\documentclass{article}
\usepackage{amsmath}
\usepackage{amsfonts}
\usepackage{amssymb}
\usepackage{amsthm}  

\newtheorem{theorem}{Theorem}
\newtheorem{proposition}{Proposition}
\newtheorem{conjecture}{Conjecture}

\title{Polynomial Potential Minimization on the Unit Circle and Related Optimization Problems}
\author{Josiah Micah Park}
\date{}

\begin{document}

\maketitle

\section{Introduction}

In the following, we study the minimization of polynomial potentials \( f(t) \) on the unit circle, where the potentials take the form 
\[
f(t) = \sum_{i=1}^n b_i x^{2i}, \quad b_i \in \mathbb{R}.
\]
This form arises in the context of truncations of expansions of \( p \)-frame potentials. One approach to minimize these potentials involves rewriting the integral as a sum of integrals obtained by expanding the potential \( f(t) = \sum_{i=1}^n c_i T_i(t) \) in terms of Chebyshev polynomials. By replacing the inner product \( \langle x, y \rangle \) with \( \cos(\theta_{x, y}) \), we can reformulate the original problem as:

\[
\min_{\mu \in P(T)} \int_T \int_T f(\langle x, y \rangle) d\mu(x) d\mu(y)
\]
as an equivalent form:
\[
\min_{\nu \in P([-\pi, \pi])} \sum_{i=1}^n c_i \int_{-\pi}^\pi \int_{-\pi}^\pi \cos(n(x - y)) d\nu(x) d\nu(y).
\]
By using the addition formula for trigonometric functions, this becomes:
\[
\min_{\nu \in P([-\pi, \pi])} \sum_{i=1}^n c_i \sum_k a_k \left[ \int_{-\pi}^\pi g_k(t) d\nu(t) \right]^2,
\]
where \( g_k(t) \) are functions selected from \( \sin(k t) \) and \( \cos(k t) \). This formulation is based on repeated use of the addition formula in trigonometric identities.

\section{Classical Results on Probability Measures and Toeplitz Matrices}

Probability measures on the unit circle are in one-to-one correspondence with certain Toeplitz matrices, as stated in the following classical result from Shohat, Tamarkin, 1937.

\begin{theorem}
A necessary condition for the solution of the trigonometric moment problem
\[
\mu_n = \int_0^{2\pi} e^{in\theta} d\mu, \quad n = 0, \pm 1, \pm 2, \dots, \quad \mu_{-n} = \mu_n
\]
is that all Toeplitz forms
\[
\sum_{j,l=0}^{n} \mu_{j-l} x_j x_l \geq 0, \quad n = 0, 1, 2, \dots
\]
are satisfied.
\end{theorem}

This result allows us to rewrite the problem as:
\[
P^* = \min_{\nu} \sum_i P(\nu_i, a_i, c_i) \quad \text{subject to} \quad A\nu \succeq 0,
\]
where the matrix \( A \) is given by:
\[
A\nu = \begin{pmatrix}
1 & \nu_1 & \nu_2 & \dots & \nu_n \\
\nu_{-1} & 1 & \nu_1 & \dots & \nu_{n-1} \\
\nu_{-2} & \nu_{-1} & 1 & \dots & \nu_{n-2} \\
\vdots & \vdots & \vdots & \ddots & \vdots \\
\nu_{-n} & \nu_{-(n-1)} & \nu_{-(n-2)} & \dots & 1
\end{pmatrix}.
\]
The function being minimized is a polynomial in the \( \nu_i \), while the constraints are semialgebraic. This suggests that for rational coefficients \( a_i \) and \( c_i \), the optimal value is algebraic. 

\begin{proposition}
For rational coefficients, \( P^* \), the optimal value of the polynomial potential minimization problem, is algebraic.
\end{proposition}

\begin{proof}
Let \( S = \left\{ \{\nu_i\}_{i=-n}^{n}, A\nu \succeq 0, \nu_0 = 1 \right\} \) denote the first moments for a probability measure supported on the circle. Since \( S \) is a semi-algebraic set, so is \( T = \left\{ (\nu, x) \mid \nu \in S, x \geq P_{a, c}(\nu) \right\} \). By the Tarski-Seidenberg principle, the projection of a semi-algebraic set is also semi-algebraic. Thus, \( P^* = \min(\text{proj}_x T) \) is algebraic.
\end{proof}

\section{Energy Minimization on the Sphere}

Next, we examine the energy minimization problem for polynomial potentials \( F(t) = \sum_{k=0}^m \alpha_k t^k \) on spheres. Let \( S = \{Q \in \mathbb{R}^{N \times N} \mid Q \succeq 0, \text{diag}(Q) = 1, \text{rank}(Q) \leq d + 1 \} \) denote a semi-algebraic set of matrices. Measures \( \nu \) of support at most \( N \) with equal masses correspond to matrices \( Q \in S \) (by abuse of notation). The energy integral for \( \nu \) is given by:
\[
I_F(\nu) = \frac{1}{N^2} \sum_{i,j=1}^N F(Q_{i,j}).
\]
Thus, the energy minimization problem for \( N \) equally weighted points on the sphere \( S^d \) is algebraic, and the optimal value \( Q^* \) is algebraic.

\begin{proposition}
For polynomial potentials \( F(t) = \sum_{k=0}^m \alpha_k t^k \), where \( \alpha_k \) are rational, the optimal value \( Q^* \) for the energy minimization problem on the sphere \( S^d \) is algebraic.
\end{proposition}

\begin{proof}
Since \( S \) is a semi-algebraic set and the set 
$$ T = \left\{ (Q, x) \mid Q \in S, \nu = (Q, \omega), x \geq I_F(\nu) = \frac{1}{N^2} \sum_{i,j} \alpha_k Q_{i,j}^k \right\} $$
is also semi-algebraic, the projection \( \text{proj}_x T \) is semi-algebraic by the Tarski-Seidenberg principle. Therefore, \( \min(\text{proj}_x T) \) is algebraic, implying that the optimal value \( Q^* \) is algebraic.
\end{proof}

\section{Minimization of Potential \( t^{2k} + \alpha C_{2k}^2(t) \) and \( p \)-Frame Potential}

This section examines the minimization of two potentials: \( t^{2k} + \alpha C_{2k}^2(t) \) and the \( p \)-frame potential \( |t|^p \). Numerical evidence suggests that minimizers of the first potential as \( \alpha \to 0 \) are related to minimizers of the second potential as \( p \to 2k \). More precisely, it is conjectured that:
\[
\mu_{\alpha \to 0} = \lim_{p \to 2k} \arg \min_{\mu} \int \int | \langle x, y \rangle |^p d\mu(x) d\mu(y),
\]
and that the measures \( \mu_{\alpha \to 0} \) and \( \nu_{p \to 2k} \) are related.

We expect that both problems have solutions with small weighted designs, supported on discrete sets of points, as \( \alpha \to 0 \) and \( p \to 2k \). The behavior of the coefficients in the expansion of the \( p \)-frame potential into Gegenbauer polynomials also shows that the minimizers of these potentials are closely related.

\section{Minimization of Polynomial Potentials on the Unit Circle}

For the two-dimensional case on the unit circle, equally distributed measures with even support sizes have Chebyshev polynomial averages. If the set of measures is constrained to ones of this form, minimizing the potentials \( f(t) = \sum_j \alpha_j C_j(t) + \sum_k \beta_k C_k(t) \) corresponds to finding the smallest sub-sum (over multiples of \( 2m \)) in the sequence of coefficients given by the potential.

This minimization problem does not always have unique solutions. For example, the function 
\[
f(t) = T_0(t) + T_2(t) + T_4(t) + 3T_6(t)
\]
attains its minimum for two distinct measures: one with support on antipodal points and another with support on a regular hexagon. However, numerically, a more optimal configuration is found for a measure supported on two points with an inner product close to a root of the derivative of \( f(t) \).

\section{Conjectures and Future Work}

We make a conjecture based on the behavior of the minimizers for these polynomial potentials:

\begin{conjecture}
If \( f(t) = \sum_j \alpha_j C_j(t) + \sum_k \beta_k C_k(t) \) does not have a unique minimal sub-sum, the minimizer will be a measure supported on a set of points whose inner products correspond to the zeros of the derivative of \( f(t) \).
\end{conjecture}

\end{document}